\documentclass[12pt]{article}

\usepackage[top=1in, bottom=1in, left=1in, right=1in]{geometry}

\usepackage{amsmath,amssymb,amsthm}
\usepackage{color}

\newcommand{\R}{\mathbb{R}}
\newcommand{\Z}{\mathbb{Z}}

\newcommand{\im}{\mathrm{i}}

\newcommand{\In}{\mathrm{in}}
\newcommand{\Sc}{\mathrm{sc}}
\newcommand{\loc}{\mathrm{loc}}

\newtheorem{theorem}{Theorem}[section]

\newtheorem{corollary}{Corollary}[section]
\newtheorem*{example}{Example}
\theoremstyle{definition}
\newtheorem{remark}{Remark}[section]

\newcommand{\hh}[1]{{\color{blue} #1}}

\newcommand{\pare}[1]{\left(#1\right)}

\newcommand{\AND}{\qquad\text{and}\qquad}

\begin{document}

\title{
	{On the Regularity of Non-Scattering Anisotropic Inhomogeneities}
}

\author{Fioralba Cakoni\footnote{Department of Mathematics, Rutgers University, New Brunswick, New Jersey, USA (fc292@math.rutgers.edu)} \; Michael S. Vogelius\footnote{Department of Mathematics, Rutgers University, New Brunswick, New Jersey, USA (vogelius@math.rutgers.edu)}\; and \; Jingni Xiao \footnote{Department of Mathematics, Drexel University, Philadelphia, Pennsylvania, USA (xiaojn@live.com)}}
 
 %\date{\empty}
\maketitle
\vskip 20pt

\date{\empty}
\maketitle
  
 \begin{abstract}
 \noindent
In this paper we examine necessary conditions for an anisotropic inhomogeneous medium  to be non-scattering at a single wave number and for a single  incident field. These conditions are expressed in terms of the regularity of the boundary of the inhomogeneity. We assume that the coefficients, characterizing the constitutive material properties of the medium, are sufficiently smooth, and
the incident wave is appropriately non-degenerate. Our analysis utilizes the Hodograph transform
as well as regularity results for  
nonlinear elliptic partial differential equations. Our approach requires that the boundary a-priori is of class $C^{1,\alpha}$ for some $0<\alpha<1$.
\end{abstract} 
 
\section{Introduction}
Spectral problems with the wave number as spectral parameter play a central role in mathematical scattering theory. A particularly noteworthy development has been the theory of scattering resonances, which correspond  to the poles of  the scattering operator viewed as meromorphic operator valued function of the wave number \cite{zow2}. Quite related is the question of injectivity of the relative scattering 
{(incident-to-outgoing \cite{CK19,res1}) operator} 
for inhomogeneous media, which has led to the study of the transmission eigenvalues \cite{CCHbook}. To be precise,  given an inhomogeneous medium of compact support, the question is if there are wave numbers for which there {exist} incident waves that are not scattered by the medium, i.e.,  the medium is rendered invisible to this particular probing experiment.  Such wave numbers are referred to as {\it{non-scattering wave numbers}}, and the corresponding incident waves as non-scattering incident fields. An inhomogeneous medium  that admits non-scattering wave numbers is called {\it non-scattering inhomogeneity}.

\noindent
If the given inhomogeneity 
{does not} 
scatter an incident wave at a fixed wave number, then inside the support of the media this  incident wave and the corresponding total field solve the Helmholtz equation governing wave propagation in the background and in the medium, respectively, and share the same Cauchy data on the boundary. This homogeneous boundary value problem is known as the transmission eigenvalue problem, and its eigenvalues (the wave numbers) are referred to as {\it{transmission eigenvalues}}. The transmission eigenvalue problem is non-selfadjoint and has an interesting mathematical structure. Its  formulation involves the governing equations for the  medium and  background only inside the support of the inhomogeneity, i.e., it does not depend on the external probing wave. Transmission eigenvalues play a central role in  inverse scattering theory for inhomogeneous media; we refer the reader to \cite{CCHbook} for an up-to-date discussion of transmission eigenvalues, their applications to inverse scattering, and the vast literature on the subject.  In particular, non-scattering wave numbers form a subset (possibly empty) of the  real transmission eigenvalues. The existence of  infinitely many real transmission eigenvalues is proven for a large class of (not necessarily regular) inhomogeneities. The question then arises, whether these inhomogeneities are non-scattering, in other words whether any of the corresponding real transmission eigenvalues is a non-scattering wave number. It turns out that the existence of a non-scattering wave number (unlike the existence of real transmission eigenvalues) generically implies a certain regularity of the inhomogeneity.  Such regularity results  are proven in \cite{CakVog} and \cite{SS} for the case of  isotropic media with contrast only in the lower order term of the Helmholtz equation.   In this paper we  deal with scattering by an anisotropic medium governed by the Helmholtz equation, with variable coefficients in the principal operator and the lower order term. We establish necessary conditions for such a medium  to be non-scattering, or equivalently by negation, sufficient conditions for it to be scattering. The investigation of non-scattering anisotropic media presents additional mathematical difficulties and leads to interesting open questions.

\bigskip

\noindent
We consider the following time-harmonic scattering problem 
\begin{equation}\label{MainGov1}
	\begin{split}
		\nabla\cdot A \nabla u + k^2 n  u=0&\quad\mbox{in $~\R^d$}, 
		\\u=u^{\In}+u^{\Sc}&\quad\mbox{in $~\R^d$},
		\\\lim_{|x|\to\infty} |x|^{\frac{d-1}{2}}\left(\frac{\partial }{\partial |x|}u^{\Sc}-\im ku^{\Sc}\right)=0, &~\mbox{ uniformly for all $\hat{x}:=\frac{x}{|x|}\in \mathcal{S}^{d-1}$},
	\end{split}
\end{equation}
where $k>0$ is the wave number proportional to interrogating frequencies, and $d\geq 2$. Here, $u^{\In}$ is a given incident field which satisfies the Helmholtz equation
$\Delta u^{\In}+k^2u^{\In}=0$ in $\R^d$ and $u^{\Sc}$ is the corresponding scattered field. We assume that  $n$ is a scalar function in $L^{\infty}(\R^d)$, and $A=(a_{ij})$ is an $d\times d$ symmetric matrix-valued function  with $L^{\infty}(\R^d)$ entries, satisfying 
\begin{equation}\label{eq:Aellip}
c_0^{-1}|\xi|^2 \le\xi^{\top} A(x) \xi\le c_0\,|\xi|^2\qquad\mbox{for almost  all $x\in\mathbb{R}^d$ and all $\xi\in\mathbb{R}^d$},
\end{equation}
for some positive constant $c_0$. We further assume that $A-I$ and $n-1$ are supported in $\overline{\Omega} \subset {\mathbb R}^d$, where $\Omega$  is a bounded Lipschitz domain.  Let $\nu$ denote the outward unit normal vector to $\partial \Omega$. In this model $A(x),n(x)$, for $x\in \Omega$, characterize the constitutive material properties of an anisotropic dielectric inhomogeneous medium occupying the region $\Omega$ and located in an isotropic homogeneous background with constitutive material properties scaled to one. In what follows $(A,n,\Omega)$ refers to this inhomogeneous medium. It is known that  the scattering problem (\ref{MainGov1}) admits a \emph{unique solution} $u\in H^1_{\loc}(\R^d)$ (see e.g. \cite{CCHbook}).  Note: the fact that $u\in H^1_{\loc}(\R^d)$ solves (\ref{MainGov1}) implies 
\begin{equation}
u^-=u^+,\quad   \qquad \nu^{\top} A \nabla u^-=\partial_{\nu}u^+,\qquad\quad \mbox{on $\partial\Omega$}~,
\end{equation}
where $+$ and $-$ indicate traces on the boundary $\partial \Omega$ from outside and inside $\Omega$, respectively.  Additional regularity conditions on the inhomogeneity $(A,n,\Omega)$ will be imposed later as required. Of main interest to us is the case when $A$ has a jump across $\partial \Omega$. In the following, when $A$ (or $n$) is continuous on $\overline{\Omega}$, and we talk about the value of $A$ (or $n$) on $\partial\Omega$, we always mean the limit value from inside $\Omega$.

\noindent
The main focus of this paper is to investigate the implied regularity on $\partial \Omega$ if it happens that  $u^{\In}$ is a non-scattering wave for the medium $(A,n, \Omega)$, in other words if  $u^{\Sc}$ vanishes identically outside $\Omega$. In this case, with $v:=u^{\In}$, for simplicity of notation, we have that $u|_{\Omega}$ along with $v$ becomes a solution to the following homogeneous problem
\begin{equation}\label{ITEP}
	\begin{split}
		\nabla\cdot A \nabla u + k^2 n  u=0,\quad  \Delta v + k^2 v=0,&\qquad\mbox{in $\Omega$}~,\\
		u=v,\quad   \nu^{\top} A \nabla u=\partial_{\nu}v,&\qquad\mbox{on $\partial\Omega$}~,
	\end{split}
\end{equation}
with $v$ actually being a solution to 
\begin{equation}\label{eqv}
\mbox{$ \Delta v + k^2 v=0$ \; in all of $\R^d$}.
\end{equation}
We note that $v$ is real analytic in $\R^d$. In our earlier terminology, the wave number $k>0$ is a non-scattering wave number for $(A,n, \Omega)$.  Considering only the set of homogeneous equations  (\ref{ITEP}), also known as the transmission eigenvalue problem, we conclude that  $k>0$ is a {\it real} transmission eigenvalue. The anisotropic transmission eigenvalue problem has been subject of extensive investigation  in the past decade  \cite{CCHbook}. Despite its deceptively simple formulation, it is a  non-selfadjoint  eigenvalue problem even for real valued  coefficients $A$ and $n$. Nevertheless, for $A,n$ real valued, and in $L^\infty(\Omega)$ and $\partial \Omega$ Lipschitz, the existence of an infinite discrete set of real transmission eigenvalues accumulating only at $+\infty$ is proven in \cite{CGH} provided $A-I$ is one-sign definite uniformly  in $\Omega$ and $n\equiv 1$, and in \cite{CK} provided both $A-I$ and $n-1$ are one-sign definite  (the same or opposite sign)  uniformly  in $\Omega$. The state-of the-art of the spectral analysis for the transmission eigenvalue problem, including discreteness of real and complex eigenvalues, completeness of generalized eigenfunctions in $(L^2(\Omega))^2$, and Weyl asymptotics for the eigenvalue counting function, can be found in \cite{nguyen}. This spectral analysis is conditional on some 
{``ellipticity''} 
assumptions on the coefficients at  the boundary $\partial \Omega$. More specifically, it requires that $\partial \Omega$ is of class $C^2$, and that $A$ and $n$ are continuous on $\overline\Omega$ and satisfy  for all $x\in \partial \Omega$\hh{,}
$$(A(x)\nu \cdot \nu)(A(x)\xi \cdot \xi)-(A(x)\nu \cdot \xi)^2\neq
 1 \qquad \mbox{and} \qquad (A(x)\nu \cdot  \nu)n(x)\neq 1,$$
 for all unit vectors $\xi\in {\mathbb R}^d$  perpendicular to the normal vector $\nu$ (the first condition is equivalent to the complementing condition, due to Agmon, Douglis and Nirenberg \cite{ADN}).  However, under these above less restrictive assumptions on $A$ and $n$, it is not known whether real transmission eigenvalues exist. We refer the reader to \cite{vodev} for results about the location  of transmission eigenvalues in the complex plane.

\noindent
Whereas real transmission eigenvalues 
{exist} 
for a broad class of (not necessarily smooth) anisotropic inhomogeneities, the main result  of this paper states  that the existence of a non-scattering wave number associated with an appropriate incident field, for regular $A$ and $n$ implies a certain regularity of  the  inhomogeneity. In particular we prove that  if $A$ and $n$ are real analytic in $\overline\Omega$, $A\neq I$ on $\partial \Omega$, and $\partial \Omega$ is a-priori of class $C^{1,\alpha}$ for some $0<\alpha<1$, and if $\partial \Omega$ is not analytic in any neighborhood of a point $P\in \partial \Omega$,  then $(A,n,\Omega)$ scatters all incident fields having a  non-vanishing conormal derivative at that point 
{(see \eqref{eq:nondeg})}. In terms of the transmission eigenvalue problem, formulated only in $\Omega$, our result provides a necessary condition, that the part $v$ of the corresponding eigenfunction can be extended as a solution of the Helmholtz equation in the exterior of $\Omega$ ({\it i.e.}, $v$ is real analytic up to the boundary). We establish similar results for less regular $A$ and $n$. To prove these results we employ the Hodograph transform to locally straighten the boundary  and transform (\ref{ITEP}) to a strongly elliptic second order nonlinear partial differential equation in divergence form\hh{,} accompanied with a nonlinear oblique derivative boundary condition. The regularity result is then obtained by appealing to the work by Agmon, Douglis and Nirenberg \cite{ADN}. The idea to use  the Hodograph transform for $C^{1,\alpha}$  boundaries is inspired by the work of Kinderlehrer and Nirenberg \cite{pom000}  (see also Alessandrini and Isakov \cite{AlesIsa96}). Scattering from inhomogeneities $(A,n,\Omega)$ where $\partial \Omega$ contains corner singularities is investigated in \cite{CaX} and \cite{Xiao} by a different approach based on 
{the so-called%rapidly decaying 
}
Complex Geometric Optics (CGO) solutions combined with asymptotic analysis in a neighborhood of the corner. For 
the readers' convenience, we summarize these results on corner scattering in Section \ref{conc}, where we also provide a simple example of a  non-scattering inhomogeneity of the form $(A,n,\Omega)$ with corners. If $\Omega$ is a ball of radius $R$ centered at the origin, and $A:=a(r)I$, $n:=n(r)$ with scalar functions  depending only on the radial variable $r$, satisfying
$$\frac{1}{R}\int_0^R\left(\frac{n(r)}{a(r)}\right)^{1/2} \,dr \neq 1~,$$
then it is possible to show by separation of variables the existence of infinitely many non-scattering wave numbers \cite{CK, CK19}.  In fact, for this spherically stratified medium the set of transmission eigenvalues and the set of non-scattering wave numbers coincide. Furthermore, the non-scattering incident waves are superpositions of plane waves, otherwise known as Herglotz wave functions with particular densities, and each density is associated with an infinite set of non-scattering wave numbers.

\noindent
The fact that we in this paper assume that $A-I\neq 0$ on the boundary is essential, and it makes the analysis more challenging. Regularity results for  non-scattering inhomogeneities with $A\equiv I$ in ${\mathbb R}^d$ are established in \cite{CakVog} and \cite{SS} (see also \cite{liu}). In this case the starting regularity  of the boundary is Lipschitz, and the incident wave can not vanish on boundary points of interest. These results can be extended to the case when $A(x)=a(x)I$ with a scalar function  $a(x)$ which is at least $C^{2,\alpha}$ up to and across the boundary, by use of a standard Liouville transformation. However, they do not extend to the case when $a$ has a jump across the boundary, nor to the anisotropic case.
Scattering from corners has been investigated  in \cite{nonscat4,CaX,ElH18,HSV16,PSV17,Xiao}. 
It turns out that for geometries with corners and smooth $A$ across the boundary (in particular $A=I$ on $\partial \Omega$) the non-vanishing  of the incident wave is not necessary to achieve non-zero scattered field. These results provide the  foundation for  proving that  a convex polygonal inhomogeneity with smooth  $A$  in a neighborhood of all corners is uniquely determined from scattering data corresponding to one single incident wave \cite{CaX,ElH18,HSV16,Xiao}. Our example in Section \ref{conc} show that such results cannot be expected in general for $(A,n,\Omega)$, when $A$ has a jump across the boundary, without additional assumptions on $A$ and $n$. Finally let us mention, that the non-scattering phenomenon of spherically  symmetric media is very unstable. More precisely, if a disk in ${\mathbb R}^2$, with $A:=I$ and $n>0$ constant, is perturbed using a very broad class of perturbations,  then there exist at most finitely many positive wave numbers for which  Herglotz wave  functions, with  smooth non-trivial densities in an appropriately compact set, that can be non-scattering (see \cite{vx,vx2}). 
\section{The Hodograph Transform and the Main Result} 
We start  by rewriting  (\ref{ITEP})-(\ref{eqv}) in a different  form more convenient for our analysis, that is, we write it in  terms of the scattered field instead of the total field. To this end, we define $w:=u-v$. Then $w$ satisfies
\begin{equation}\label{PDEw}
\nabla\cdot A \nabla w + k^2 n  w=-\nabla\cdot (A-I) \nabla v- k^2 (n-1)  v\qquad\mbox{in $\Omega$}~,
\end{equation}
with
	\begin{equation}\label{eq:BCw}
	w=0~~ \hbox{and } ~~\nu^{\top} A \nabla w=-\nu^{\top} (A-I) \nabla v\qquad\mbox{on $\partial\Omega$}~,
	\end{equation}
with $v$ being a solution to \eqref{eqv}, and
hence real analytic in $\R^d$. As it will become clear later in the paper, what matters to our results is that $v$ is sufficiently regular on $\overline{\Omega}$.

\noindent
From now on we shall assume that $\Omega$ is a $C^{1,\alpha}$ domain and that the matrix valued function  $A$ has entries in $C^{1,\alpha} (\overline{\Omega})$\footnote{By $C^{k,\alpha}(\overline \Omega)$, $k$ an integer $\ge 0$, and $0<\alpha<1$, we understand functions that may be extended as $C^{k,\alpha}$ functions in an open neighborhood of $\overline\Omega$. The analogue of Whitney's Extension Theorem \cite{Whit34} for $C^{k,\alpha}$ functions, $0<\alpha<1$, asserts that  this definition of $C^{k,\alpha}(\overline \Omega)$ amounts to requiring that all derivatives of order less than or equal to $k$ be $\alpha$-H\"older continuous in $\Omega$, and up to the boundary $\partial \Omega$, with constants that are uniformly valid in $\Omega$.}.
Due to the regularity of $v$, a standard regularity result for elliptic equations (see Corollary 8.36  in \cite{GilTru01}) gives that $w\in C^{1,\alpha}(\overline{\Omega})$.
Moreover, since $w=0$ on $\partial\Omega$, we have $\nu=\pm \nabla w/|\nabla w|$, provided $\nabla w\neq 0$. The second condition in \eqref{eq:BCw} therefore leads to
\begin{equation}\label{eq:BCwnonu}
(\nabla w)^{\top} A \nabla w+(\nabla w)^{\top} (A-I) \nabla v=0\qquad\mbox{on $\partial\Omega$ %where $\nabla w\neq 0$
}.
\end{equation}
We  are interested in the case when $A$ has a jump across $\partial \Omega$. Now let  $P$ be a (fixed) point on $\partial \Omega$ with $A(P)\ne I$. 
We also assume that 
\begin{equation}\label{eq:nondeg}
\nu^{\top}(A-I)\nabla v(P)\ne 0~.
\end{equation}

\vskip 5pt

\vskip 5pt
\noindent
Then, recalling the boundary condition \eqref{eq:BCw} we also have
\begin{equation*}
\nu^{\top}A\nabla w\ne 0
\qquad\mbox{at $P$}.
\end{equation*}
Up to a rigid change of coordinates we may assume that $\nu^{\top}A=(-c_1,0, \ldots\, ,0)$ with $c_1=|A\nu|> 0$ at $P$, and that $P=0$.
Without loss of generality we may also assume that $\nu^{\top}(P)A(P)\nabla w(P)=-c_1\partial_1w(0)<0$.
Then, by the regularity of $w$ %(and that of $A$ and $\nu$), 
we can find positive constants $r$ and $c_2$ such that
\begin{equation}\label{eq:d1w<0}
%\nabla w(x)\neq 0
%\AnD
0<c_2^{-1}<\partial_1w(x)<c_2,
\qquad\mbox{for all $x\in\overline{\Omega}\cap B_r(0)$}~,
\end{equation}
and hence 
\begin{equation}\label{eq:dv<0}
0< c_3^{-1}<-\frac{(\nabla w)^{\top} (A-I) \nabla v}{|\nabla w|}<c_3~,
\qquad\mbox{for all $x\in\partial{\Omega}\cap B_r(0)$}
\end{equation}
for some constant $c_3$.
Denoting $x=(x_1, x_2, \ldots \, ,x_d)=(x_1, x')$, we consider the mapping
	$$
H:\	x\mapsto y = (w(x), x')~,\qquad x\in\overline{\Omega}\cap B_r(0)~.
	$$
Due to the regularity of $w$, the map $H$ is in $C^{1,\alpha}(\overline{\Omega}\cap B_r(0))$. 
	%As consequence, $H$ can be extended to the whole $B_r(0)$ as a $C^{1,\alpha}$ mapping. 
	It can be calculated directly that the Jacobian of $H$ is
	\begin{equation*}
	DH(x)=\begin{bmatrix}
	\partial_{x_1}w&\nabla_{x'}^{\top}w
	\\
	0&I
	\end{bmatrix}~.
	\end{equation*}
Therefore, $H$ is invertible on $\overline{\Omega}\cap B_r(0)$, thanks to the non-vanishing of $\partial_{x_1}w$. We may conclude that $H$ is bijective from $\Omega \cap B_r(0)$ to $V^+=V\cap\{ (y_1,y')~: ~y_1>0\}$, and also from $\partial \Omega \cap B_r(0)$ to $\Sigma =V\cap \{ (y_1,y')~: ~y_1=0\}$, where $V$ is an open neighborhood of $0$. The inverse map of $H$ can be expressed as $H^{-1}(y)=(z(y),y')$ for some function $z\in C^{1,\alpha}(V^+\cup \Sigma)$. In terms of $z$ the Jacobian of $H$ can be re-written as
	\begin{equation*}
	DH(x)=\begin{bmatrix}
	\partial_{x_1}w&\nabla_{x'}^{\top}w
	\\
	0&I
	\end{bmatrix}
=\begin{bmatrix}
	1/\partial_{y_1}z&-\nabla_{y'}^{\top}z/\partial_{y_1}z
	\\
	0&I
\end{bmatrix}\circ H(x)~.
	\end{equation*}
Next, we derive equations for $z$. By a direct calculation we obtain for $y=H(x)$ that
\begin{equation*}
	\nabla_x^{\top}=\widetilde{\nabla}^{\top}z~\partial_{y_1}+(0,\nabla_{y'}^{\top})~,
	\qquad
	\mbox{where $\widetilde{\nabla}^{\top}z=\dfrac{1}{\partial_{y_1}z}(1,-\nabla_{y'}^{\top}z)$}~.
\end{equation*}
In the following, all the differential operators are with respect to $y$, unless they carry explicit subscripts of $x$.
We then have
\begin{equation*}
	(\nabla_x w)\circ H^{-1}=\widetilde{\nabla}z
	\AND
	%\end{equation*}
	%and 
	%\begin{equation*}
	(\nabla_x\cdot A\nabla_x w)\circ H^{-1}
	=
	\pare{\widetilde{\nabla}^{\top}z~\partial_{1}
		+
		(0,\nabla_{y'}^{\top})}
	(A_H\widetilde{\nabla}z)~,
\end{equation*}
with $A_H=A\circ H^{-1}$.
A direct calculation gives that
\begin{equation*}%\label{eq:idcal}
	2\,\widetilde{\nabla}^{\top}z~\partial_{y_1}
	(A_H \widetilde{\nabla}z)
	=
	\partial_{y_1}\pare{ (\widetilde{\nabla}^{\top}z)\,  A_H \widetilde{\nabla}z} 
	+
	\widetilde{\nabla}^{\top}z\, 
	\pare{\partial_{y_1}A_H}\,\widetilde{\nabla}z~,
\end{equation*}
and that 
\begin{equation*}
\partial_{y_1}A_H(y)=	\partial_{y_1}A(H^{-1}(y))=\partial_{y_1}A(z(y),y')
	=(\partial_{y_1}z)\,\frac{\partial A}{\partial x_1}(H^{-1}(y))~.
\end{equation*}
Therefore, we deduce from \eqref{PDEw} that $z(y)$ satisfies 
\begin{equation}\label{eq:PDEz}
	\sum_{j=1}^{d}\partial_{y_j}a_j(y,z,\nabla z)+a_0(y,z,\nabla z)=0
	\qquad\mbox{in $V^+$},
\end{equation}
where $a_j$, $j=0,1,\ldots,d$, are functions of $2d+1$ arguments given by
\begin{equation*}
	a_1=\frac{1}{2}(\widetilde{\nabla}^{\top}z)A_H\widetilde{\nabla}z~,
	\qquad
	a_j=(A_H\widetilde{\nabla}z)_j~,\quad  j=2,\ldots,d~,
\end{equation*}
and 
\begin{equation*}
	a_0=\frac{1}{2}\,
	(\partial_{1}z)\,
	\widetilde{\nabla}^{\top}z\pare{\partial_{x_1}A}\widetilde{\nabla}z
	+k^2ny_1
	+\nabla_x\cdot (A-I) \nabla_x v+ k^2 (n-1)  v~.
\end{equation*}
%with the notation
%\begin{equation*}
%	\widetilde{\nabla}^{\top}z=\dfrac{1}{\partial_{1}z}(1, -\nabla_{y'}^{\top}z).
%\end{equation*}
Here, in the lower order term $a_0$, with an abuse of notation, we still use $n$ to denote $n\circ H^{-1}(y)=n(z(y),y')$; the same applies to $A$ and $v$, as well as their derivatives with respect to $x$.
Similarly, from the boundary condition \eqref{eq:BCwnonu} for $w$ we obtain that
\begin{equation}\label{eq:BCz}
	b(y',z,\nabla z)
	:=
	(\widetilde{\nabla}z)^{\top}A_H\widetilde{\nabla}z
	+(\widetilde{\nabla}z)^{\top}(A_H-I)\nabla_x v
	=0\qquad\mbox{on $\Sigma$}~,
\end{equation}
where $A_H=A_H(0,y')$ and $\nabla_x v=(\nabla_x v)\circ H^{-1}(0,y')$~.

We claim that \eqref{eq:PDEz} is a uniformly strongly (nonlinear) elliptic equation and \eqref{eq:BCz} is an associated proper (nonlinear) oblique derivative boundary condition. 
This result heavily depends on the fact that $\partial_{1}z=1/\partial_{x_1}w$ and thus 
\begin{equation}\label{eq:d1z>0}
0<c_2^{-1}<\partial_1z(y)<c_2
\qquad\mbox{for all $y\in V^+\cup \Sigma$}~.
\end{equation}
We shall present the proof of this claim in the Appendix.

\begin{remark}
In the special case when $A={a}I$ in $\overline{\Omega}$ with ${a}$ being a positive constant, the problem for $z$ becomes
	\begin{equation*}
		\begin{split}
			-\frac{1}{2}\,{a}\,\partial_{1}\pare{\frac{1}{(\partial_{1}z)^2}(1+|\nabla_{y'}z|^2)}
		+{a}\nabla_{y'}\cdot\pare{ \frac{1}{\partial_{1}z}\nabla_{y'}z}
		=k^2ny_1
		+ k^2 (n-{a})  v
	&	\qquad\mbox{in $V^+$},
	\\
{a}\,\frac{1+|\nabla_{y'}z|^2}{(\partial_{1}z)^2}
+({a}-1)\frac{\partial_{x_1}v-(\nabla_{y'}z)\cdot\nabla_{x'}v}{\partial_{1}z}
=0
&\qquad\mbox{on $\Sigma$}~.
		\end{split}
	\end{equation*}
	Moreover, the non-degeneracy condition \eqref{eq:nondeg} simplifies to 
	\begin{equation*}
	a\neq 1
	\AND
	\partial_{\nu} v(P)\ne 0~.
	\end{equation*}
	
\end{remark}
\noindent
We now assume that $n\in C^{1,\alpha}(\overline{\Omega})$ and $A\in \left(C^{2,\alpha}(\overline{\Omega})\right)^{d\times d}$. 
We are then in a position to use Theorem 11.2 of \cite{ADN}, with $l=m=m_1=1$ and $p=2$,  to conclude that a solution to the boundary value problem \eqref{eq:PDEz}-\eqref{eq:BCz} which belongs to $C^{1,\alpha}(V^+\cup\Sigma)$ is indeed in $C^{2,\alpha}(V^+\cup\Sigma)$. Here we use (among other things) that the functions $a_j$ in \eqref{eq:PDEz}, $j=0,1,\ldots,d$, are $C^{1,\alpha}$ in all the arguments, and so is $b$ in \eqref{eq:BCz}. Once we have established that $z$ is $C^{2,\alpha}$ near $P=0$, it follows that $\partial \Omega$ and $w$ are $C^{2,\alpha}$ near $0$. 
Furthermore, if $A$ and $n$ admit higher regularity up to $C^{\infty}$, by a bootstrap procedure we can establish higher regularity results for $\partial\Omega$ up to $C^{\infty}$ (as well as for $w$ and $u$).
In addition, if $A$ and $n$ are analytic, we can also prove analyticity of $\partial \Omega$ near $P=0$, using results of Morrey \cite{Morrey}. 
To be quite precise we have established
\begin{theorem}
	\label{nonscat}
	Let $\Omega$ be a bounded domain in $\R^d$, $d\ge2$, of class $C^{1,\alpha}$. Suppose for some integer $\ell\ge 1$ that $n\in C^{\ell,\alpha}(\overline{\Omega})$  and that $A\in \left(C^{\ell+1,\alpha}(\overline{\Omega})\right)^{d\times d}$ with the condition \eqref{eq:Aellip} satisfied. 
	Let $(u, v)$ be a solution to \eqref{ITEP}, with $v$ a solution to $\Delta v+k^2v=0$ in a neighborhood of $\Omega$. Let $P$ be a point on $\partial \Omega$ at which the non-degeneracy condition \eqref{eq:nondeg} is satisfied.
	Then $\partial \Omega$ is of class $C^{\ell+1,\alpha}$ near $P$.
	Moreover, if both $A$ and $n$ are $C^{\infty}$ on $\overline{\Omega}$, then $\partial\Omega$ is $C^{\infty}$ near $P$. If $A$ and $n$ are real analytic on $\overline{\Omega}$, then $\partial\Omega$ is also analytic near $P$.
\end{theorem}

\noindent
\begin{remark}
We note that we never specifically use in the  proof  of Theorem \ref{nonscat} that $v$ solves the Helmholtz solution in the exterior of $\Omega$, it suffices that it has a smooth extension. In fact, we need only sufficient regularity of $v$, that is $C^{\ell+2,\alpha}$ up to the boundary.  Furthermore, the higher regularity assumptions on the coefficients in Theorem \ref{nonscat}  are only needed locally  in $\overline{\Omega}\cap B_R(P)$ for some $R>0$.
\end{remark}

\noindent
\begin{remark}
{Theorem \ref{nonscat} in particular implies $C^{\ell+1,\alpha}$ regularity of a non-scattering inhomogeneity, provided  that the associated incident field $v$ satisfies \eqref{eq:nondeg}, since in that case the incident field  $v$ is an analytic solution to $\Delta v+k^2v=0$ in all of $\R^d$.}
\end{remark}
\noindent
Of course Theorem \ref{nonscat} only add insight if the wave number $k$ is a real transmission eigenvalue (which is a necessary condition for the incident field to produce a vanishing scattered field). 
%At any  $k$ other than a transmission eigenvalue, every incident field is scattered by the given %inhomogeneity. However, it is important to emphasize that we do not need to know a priori that $k>0$ is %a transmission eigenvalue, and therefore our results hold under weaker conditions on the contrast than %those (currently) needed to prove the existence  of real transmission eigenvalues. 
\noindent
If $k>0$ is a transmission eigenvalue, a variant of Theorem \ref{nonscat} also sheds light onto regularity up to the boundary of the $v$ part of the corresponding transmission eigenfunction.   In particular we have the following result.
\begin{corollary}\label{cor1}
Assume $k>0$ is a real transmission eigenvalue with eigenfunction $(u, v)$ satisfying  \eqref{ITEP}. Assume that  $\partial \Omega$ is $C^{1,\alpha}$,  $A\in \left(C^{1,\alpha}(\overline{\Omega})\right)^{d\times d}$, $0<\alpha<1$ satisfying  \eqref{eq:Aellip}. Let $P\in \partial \Omega$, and assume that $A(P)-I\neq 0$. The following assertions hold:
\begin{enumerate}
\item  If $A$ and $n$  are  real analytic in a neighborhood of $P$,  and $\partial \Omega \cap B_r(P)$ is not real analytic for any ball $B_r(P)$, then $v$ can not be real analytic  in any neighborhood of $P$, unless $\nu^{\top}(A-I)\nabla v(P)= 0$.
\item If $n\in C^{\ell,\alpha}(\overline{\Omega})\cap B_R(P)$, $A\in C^{\ell+1,\alpha}(\overline{\Omega}\cap B_R(P))$ for some integer $\ell\ge 1$, $0<\alpha<1$ and some ball $B_R(P)$,  and  $\partial \Omega\cap B_r(P)$  is not $C^{\ell+1,\alpha}$  for any ball $B_r(P)$ then $v$ cannot be in $C^{\ell+1,\alpha}(\overline{\Omega}\cap B_r(P))$  for any ball $B_r(P)$, unless $\nu^{\top}(A-I)\nabla v(P)= 0$.
\end{enumerate}
\end{corollary}

\noindent
The idea to use the Hodograph transform, as an essential tool when starting from a $C^{1,\alpha}$ boundary, is one we have borrowed from Kinderlehrer and Nirenberg \cite{pom000} (see also, Alessandrini and Isakov \cite{AlesIsa96}). For the problem treated in  \cite{AlesIsa96} in the two dimensional case there is a direct approach using conformal mapping which works when starting only with a Lipschitz assumption on $\partial \Omega$. That approach does not work for our problem even in ${\mathbb R}^2$, and so it remains an open and very interesting problem to resolve whether the initial $C^{1,\alpha}$ assumption on $\partial \Omega$ can be relaxed to a Lipschitz assumption.

\medskip 
\noindent
	In the case when $A=I$ but $n$ differs from $1$, a result similar to that in Theorem \ref{nonscat} was established in \cite{CakVog} (see also \cite{HSV16}), starting from a Lipschitz domain, and under the assumption that the incident wave does not vanish  (locally). In that case, techniques of Caffarelli \cite{pom00} and Williams \cite{pom3} were used to proceed from Lipschitz to $C^{1,\alpha}$.  The result in \cite{CakVog} is entirely local in character, and therefore also applies if $A=I$ in a neighborhood of $P\in \partial \Omega$ with $n(P)\ne 1$. In that case it only guarantees local regularity of $\partial \Omega$ near $P$. For $A=I$ the non-vanishing of the incident wave  turns out not to be necessary for scattering from corners; see e.g.  \cite{CaX,ElH18,PSV17,Xiao}.

\section{Concluding Remarks}\label{conc}
\subsection{Examples of Non-Scattering Inhomogeneities}
Scattering properties of inhomogeneities $(A,n,\Omega)$ in ${\mathbb R}^2$ with support $\Omega$ containing a corner are studied in \cite{CaX,Xiao}.  In the remark below we summarize the main results obtained in \cite{Xiao} in the case, when $A$ is isotropic and is discontinuous across $\partial \Omega$.
\begin{remark}\label{jingni}
{Let $\Omega \in {\mathbb R}^2$ and  $P\in \partial \Omega$ be the vertex of a corner with aperture $2\theta$, $\theta\in(0,\pi/2)\cup(\pi/2,\pi)$. Assume that $A:=aI$ with a scalar function $a\in C^{1,\alpha}(\overline{\Omega})\cap B_R(P)$, for some ball $B_R(P)$, such that $a(P)\neq 1$. Then this inhomogeneity $(A,n,\Omega)$ scatters any incident field $v$ (which is a solution to the Helmholtz equation in a neighborhood of $\Omega)$, provided that:
\begin{itemize}
\item [(a)] $v(P)\neq 0$ and $\nabla v(P)\neq 0$, or
\item[(b)]  $v(P)=0$ and $2\theta\neq p\pi/N$ for any integer $p\geq 1$, with $N$ being the vanishing order of $v$ at $P$ (in particular if $2\theta\in (0,2\pi)/{\mathbb Q}\pi)$), or 
\item [(c)] $v(P)\neq 0$ and $\nabla v(P)= 0$, and  either $2\theta \notin \left\{\pi/2, 3\pi/2\right\}$ 
{with $v$ admitting a specific expansion,} 
or $2\theta \in \left\{\pi/2, 3\pi/2\right\}$ and $a(P)\neq n(P)$.
\end{itemize}
}
\end{remark} 
\noindent
In the next example we show that (exceptional) non-scattering may indeed happen for nonsmooth inhomogeneities.
\begin{example}[{Non-Scattering Waves and Non-Scattering Isotropic Inhomogeneities with Corners}]
{\em Consider the case when $A=aI$ and $n=a$, with $a$ a constant different from $1$. In that case any $k^2$ for which \eqref{ITEP} has a non-trivial solution is either a Dirichlet eigenvalue for $-\Delta$ or a Neumann eigenvalue for $-\Delta$ on $\Omega$  \cite{CHX}. To see this, define $w_1=u-v$ and $w_2=au-v$. Then $\Delta w_j+k^2w_j=0$ in $\Omega$, $j=1,2$, and $w_1=0$ on $\partial \Omega$ while $\partial_\nu w_2=0$ on $\partial \Omega$. The functions $w_1$ and $w_2$ cannot both vanish identically, since $(u,v)$ is a non-trivial solution to (\ref{ITEP}), and so it follows that $k^2$ is {\it either} a Dirichlet eigenvalue {\it or} a Neumann eigenvalue for $-\Delta$. Conversely, if $\lambda =k^2$ is a Dirichlet or a Neumann eigenvalue for $-\Delta$ on $\Omega$ with associated eigenfunction $w$, then $(u,v)=(w,aw)$ or $(u,v)=(w,w)$, respectively, is a nontrivial solution to
\begin{eqnarray*}
	\nabla\cdot a \nabla u + k^2 a  u=0,\quad  \Delta v + k^2 v=0,&&\qquad\mbox{in $\Omega$}~,\\
	u=v,\quad   a \partial_{\nu} u=\partial_{\nu}v,&&\qquad\mbox{on $\partial\Omega$}~.
\end{eqnarray*}
Whether $w$ corresponds to a non-scattering incident wave at wave number $k$, depends on whether $w$ can be extended to a (smooth) solution to $\Delta w+k^2w=0$ in {\it all of } $\R^d$.

\noindent
Let us now take $\Omega \subset \R^2$ be to the unit square $\Omega =Q=(0,1)\times(0,1)$. In that case all the Dirichlet eigenvalues of $-\Delta$ are given by $(p^2+q^2)\pi^2$, where $p$ and $q$ are positive integers; corresponding Dirichlet eigenfunctions are given by $w(x_1,x_2)=\sin p \pi x_1 \, \sin q \pi x_2$. All the Neumann eigenvalues of $-\Delta$ are given by $(p^2+q^2)\pi^2$, where $p$ and $q$ are non-negative integers, and corresponding Neumann eigenfunctions are given by $w(x_1,x_2) =\cos p\pi x_1\,\cos q \pi x_2$. Notice that all these eigenfunctions are smooth functions on all of $\R^2$. As a consequence it follows that, as an incident wave for the medium $0<a=n\neq 1$, $u^{\In}:=v=\cos p \pi x_1 \, \cos q \pi x_2$ is non-scattering  (at wave number $k= \sqrt{(p^2+q^2)\pi^2}$)) for any $p,q \in \Z$. The same holds for the incident wave $\sin p \pi x_1 \, \sin q \pi x_2$ for any $p,q \in \Z \setminus \{0\}$. We note that $\nabla w$ vanishes at all the corners of $\Omega$ for any of these $w$, which violates the non-degeneracy condition \eqref{eq:nondeg}. 
This example of non-scattering waves in fact falls into the possible set of non-scattering waves in Remark \ref{jingni}}. 
\end{example}
\vskip 5pt

\noindent
This very simple example illustrates two points: (i) It furnishes an example of a Lipschitz domain and isotropic $A$ for which there exist plenty of non-scattering incident waves, however, notice that these non-scattering incident waves all have vanishing gradients at the irregular boundary points. (ii) In doing so, it indicates the need for a non-degeneracy condition of the form \eqref{eq:nondeg} (or $\nabla v \ne 0$ at points, where the normal is not well defined) if one wants to prove regularity of $\partial \Omega$, starting with a Lipschitz assumption.

\vskip 5pt

\begin{example}[{Non-Scattering Inhomogeneities with Anisotropic $A$}]
{\em If one allows anisotropic $A$, there is a natural way to construct a media which is non-scattering for any incident wave. Simply let $\Phi$ be a sufficiently smooth ($C^{3,\alpha}$) diffeomorphism from $\Omega$ onto $\Omega$, with $\Phi(x)=x$ on $\partial \Omega$, and define
$$
A = \Phi_*I = \frac{D\Phi D\Phi^{\top}}{|\det D \Phi |} \circ \Phi^{-1} ~~\hbox{ and }~~ n =\Phi_*1 = \frac{1}{|\det D\Phi|}\circ \Phi^{-1}~.
$$
$(A,n.\Omega)$ is the so-called pushforward of $(I,1,\Omega)$ under $\Phi$. It is well known that the operator $\nabla \cdot (A \nabla \cdot) + k^2 n \cdot$ has the same Dirichlet-to-Neumann data map (or Cauchy data) as $\Delta + k^2 n$ for any $k>0$ (see \cite{KV,KOVW}). Therefore $(A,n,\Omega)$ is non-scattering for any wave number and any incident wave. We note that $\Phi_*I$ is anisotropic, unless $\Phi=id$, in which case $A=I$ and $n=1$. To understand why this construction does not contradict Theorem \ref{nonscat}, one must understand that if $\partial \Omega$ is not of class $C^{l+1,\alpha}$ near $P$, then either $A=\Phi_*I (P)$ is the identity matrix (remember, the gradient field of the incident wave at $P$ can be arbitrary) or $(A,n)=(\Phi_*I,\Phi_*1)$ fails to be $(C^{l+1,\alpha}, C^{l,\alpha})$ near $P$. To further illustrate this: if one takes $\Omega \subset \R^2$ be to be a polygon, with vertices $P_i$, $i=1,\ldots , N,$ and if $\Psi(x)$ is a nontrivial $C^{\infty}$ vector field that vanishes on $\partial \Omega$, then $\Phi(x)=Ix+\epsilon\Psi(x)$ is a $C^\infty$ diffeomorphism of $\Omega$ onto $\Omega$ for $\epsilon$ sufficiently small, with $\Phi(x)=x$ on $\partial \Omega$ and  with $D\Phi (P_i)=I$ for $i=1,\ldots,N$. The corresponding medium $(A,n,\Omega)=(\Phi_*I,\Phi_*1,\Omega)$ is anisotropic and non-scattering for any wave number and any incident wave; it has $C^\infty$ coefficients $A$ and $n$, however, it has $A(P_i)=I$, $n(P_i)=1$ for all $i$.}  
\end{example}

\subsection{Incident Fields satisfying the Non-Degeneracy Condition (\ref{eq:nondeg}) }\label{deg}
A natural question to ask is whether there are incident fields that satisfy the non-degeneracy condition (\ref{eq:nondeg}) at any point on $\partial \Omega$, i.e.,  whether one can find $v$ that satisfies 
\begin{equation}\label{nonA}
\Delta v+k^2v=0 \quad \mbox{in} \; {\mathbb R}^d\qquad \mbox{and} \qquad \nu^{\top}(A-I)\nabla v\ne 0 \quad \mbox{on} \; \partial \Omega~,
\end{equation}
provided that $A\neq I$ on $\partial \Omega$.  For (isotropic) inhomogeneities with  $A=aI$ and $a(x)\neq 1$ for $x\in \partial \Omega$  this reduces to finding $v$ that solves
\begin{equation}\label{non}
\Delta v+k^2v=0 \quad \mbox{in} \; {\mathbb R}^d\qquad \mbox{and} \qquad \frac{\partial v}{\partial \nu}\ne 0 \quad \mbox{on} \; \partial \Omega~,
\end{equation}
A plane wave $v(x)=e^{ik\xi\cdot x}$, where $\xi\in {\mathcal S}^{d-1}$ is the unit propagation direction, is one of the most commonly used incident waves in scattering theory. Although it has non-vanishing gradient $ik\xi e^{ik\xi\cdot x}$ it doesn't satisfy the non-vanishing condition (\ref{nonA}) on all of $\partial \Omega$. In fact, for a smooth boundary, by continuity of  $\nu(x)$ there always exists a point $P\in \partial \Omega$ where $\nu^\top(P)\cdot \xi= 0$ and hence $\partial v/\partial \nu (P)=0$. This applies to point sources as well. However, there are plenty of Herglotz wave functions (superpositions of plane waves) given by
\begin{equation}\label{eqh}
v_g(x):=\int_{{\mathcal S}^{d-1}}g(\xi)e^{ik\xi\cdot x}\,ds_\xi
\end{equation}
 which satisfy (\ref{nonA}), and thus  (\ref{non}) for $A=aI$, if we exclude some values of wave number $k>0$.  Such values of $k>0$ correspond to $k^2$ being an eigenvalue of the oblique Neumann eigenvalue problem 
$$\Delta q+k^2q=0 \quad \mbox{in} \; \Omega \qquad \mbox{and} \qquad \nu^{\top}(A-I)\nabla q= 0 \quad \mbox{on} \; \partial \Omega~.$$
If we assume the regular oblique derivative condition, i.e., 
$$(A(x)\nu \cdot \nu)(A(x)\xi \cdot \xi)-(A(x)\nu \cdot \xi)^2\neq
 1 \qquad \mbox{for all} \quad \xi\in {\mathcal S}^{d-1} \quad \mbox{such that}\; \nu\cdot \xi=0~,$$
then the oblique Neumann eigenvalues form a discrete set that accumulate at $\infty$  (see e.g. \cite{GilTru01}). Now if $k^2$ is not an oblique Neumann eigenvalue, then one can find a $v\in C^{1,\alpha}(\overline{\Omega})$ that satisfies the Helmholtz equation $\Delta v+k^2v=0$ in $\Omega$ such that $\nu^{\top}(A-I)\nabla v=\gamma>0$ on $\partial \Omega$. Therefore any Herglotz wave function (\ref{eqh}) that approximates  $v$ sufficiently well in the $C^1(\overline{\Omega})$ norm\footnote{It is possible to slightly modify the duality argument in the proof of Preposition 3.4 in \cite{SS}  to show that Herglotz wave functions  (\ref{eqh}) with $g\in C^\infty({\mathcal S}^{d-1})$ are dense in the space is dense in the space $\left\{v:\in W^{1+\sigma,p}(\Omega): \, \Delta v +k^2 v=0\right\}$ for some $0<\sigma<1$ and any $p>1$, with respect to the $W^{1+\sigma,p}(\Omega)$-norm. Then the approximation property in $C^1(\overline\Omega)$ follows from the Sobolev Imbedding Theorem.} satisfies the non-degeneracy condition (\ref{nonA}) on all 
of $\partial \Omega$. Since Theorem \ref{nonscat} is meaningful only  if $k$ is a real transmission eigenvalue, the excluded values are real transmission eigenvalues $k>0$  
such that $k^2$  is simultaneously a (oblique) Neumann eigenvalue. It is an open problem if such a set of wave numbers is finite or not. 
\section*{Appendix}
Here we show that the first variation of the nonlinear partial differential equations \eqref{eq:PDEz}--\eqref{eq:BCz} is uniformly strongly elliptic with a proper oblique derivative boundary condition; this  allows us to apply Theorem 11.2 in \cite{ADN}. The first variation 
$$\mathcal{L}_z(y,\phi,\nabla \phi)=0$$
of the set of nonlinear equations \eqref{eq:PDEz}--\eqref{eq:BCz}, in shorthand written as
$$\mathcal{F}(y,z,\nabla z)=0~,$$ is defined by means of  
$$\mbox{ $\mathcal{L}_z=\mathcal{I}'(0)$ \qquad where \qquad  $\mathcal{I}(\tau)=\mathcal{F}(y,z+\tau \phi,\nabla (z+\tau \phi))$}~.$$
See \cite[Chapter 8]{Evans} and  also \cite[Page 684]{ADN}. Assuming that $z$ is a function of $d$ variables,  $\nabla z$ is regarded as $d$ arguments, namely $\partial_{j}z$, $j=1,\ldots,d$, of $\mathcal{F}$.
 
\noindent
A straightforward calculation shows that the first  variation of \eqref{eq:PDEz} has the following matrix coefficient $\widetilde{A}$ in the principal (divergence form second order) part of the operator at each point $y\in V^+$,
\begin{equation*}
\widetilde{A}
=
-\frac{1}{\partial_{1} z}	
\begin{bmatrix}
(\widetilde{\nabla}^{\top}z)\,A\,\widetilde{\nabla}z
&\frac{1}{\partial_{1} z}(\vec{a}'-A_{d-1}\nabla_{y'}z)^{\top}
\\
\frac{1}{\partial_{1} z}(\vec{a}'-A_{d-1}\nabla_{y'}z)& A_{d-1}
\end{bmatrix},
\end{equation*} 
where $\vec{a}'^{\top}=(a_{12},a_{13},\ldots,a_{1d})$ and $A_{d-1}=(a_{ij})_{i,j=2}^d$, that is,
\begin{equation*}
A=\begin{bmatrix}
a_{11}&\vec{a}'^{\top}\\\vec{a}'&A_{d-1}
\end{bmatrix}~.
\end{equation*}
For simplicity of notations, we have again used $A=(a_{ij})$ to denote $A_H = A\circ H^{-1}=(a_{ij}\circ H^{-1})$.
We calculate for $\xi=(\xi_1,\xi'^{\top})^{\top}\in\mathbb{R}^d$ that%\footnote{We are really considering $\xi=(\xi_1,\xi'^{\top})^{\top}$ as a column vector.}
\begin{equation*}
-(\partial_{1} z)\,	\xi^{\top} \widetilde{A}\,\xi
=
\xi_1^2~(\widetilde{\nabla}^{\top}z)\,A\,\widetilde{\nabla}z
+2\,\xi_1\frac{1}{\partial_{1} z}(\vec{a}'-A_{n-1}\nabla_{y'}z)^{\top}\xi'
+\xi'^{\top} {A}_{n-1}\xi'
=
\widetilde{\xi}_z^{\top}A\,\widetilde{\xi}_z~,
\end{equation*}
with
\begin{equation*}
\widetilde{\xi}_z^{\top}
=
\pare{\frac{1}{\partial_{1} z}\xi_1~,~ (\xi'-\frac{\nabla_{y'}z}{\partial_{1} z}\,\xi_1)^{\top}}
=(0,\xi'^{\top})+\xi_1\widetilde{\nabla}^{\top}z~.
\end{equation*}
It then follows by continuity (and compactness) that, for all $y\in V^+$,
$$
(c_4)^{-1}|\xi|^2~\le~ -\xi^{\top} \widetilde{A}\,\xi ~\le~ c_4 |\xi|^2,
$$
with some positive constants $c_4>0$ depending on $c_0$  in \eqref{eq:Aellip}, $c_2$  in \eqref{eq:d1w<0}, and $\|\widetilde{\nabla}z\|_{V^+\cap \Sigma}=\|\nabla_x w\|_{\overline{\Omega}\cap B_r(0)}$. This verifies the uniformly strong ellipticity condition.  

\noindent
In regards to the boundary condition \eqref{eq:BCz}, the principal part of the first variation at each boundary point $y=(0,y')$ is give by $\sum_{j=1}^{d}b_j\partial_{j}$ with 

\begin{equation*}
b_1
=
-\frac{1}{\partial_{1}z}
\pare{2\,(\widetilde{\nabla}^{\top}z)\,A\,\widetilde{\nabla}z
+
(\widetilde{\nabla}^{\top}z)\,(A-I)\,\nabla_x v}
=-\frac{1}{\partial_{1}z}(\widetilde{\nabla}^{\top}z)\,A\,\widetilde{\nabla}z~,
\end{equation*}
and 
\begin{equation*}
b_j
=
-\frac{1}{\partial_{1}z}
\pare{2A\,\widetilde{\nabla}z
	+
	(A-I)\,\nabla_x v}_j,
\qquad j=2,\ldots,d~,
\end{equation*}
where again, by abuse of notation,  we use $A$ in place of $A\circ H^{-1}(0,y')$ and $\nabla_x v$ in place of $(\nabla_x v)\circ H^{-1}(0,y')$.
From \eqref{eq:Aellip} and \eqref{eq:d1z>0} we get that 
\begin{equation*}
0< c_0^{-1}c_2^{-3}<c_0^{-1}c_2^{-1}|\widetilde{\nabla}z|^2< -b_1 \qquad
\mbox{for all $y=(0,y')$ on $\Sigma$}~.
\end{equation*}
This ensures that the linearized boundary condition (on $\Sigma$) is a proper oblique derivative condition, and as a consequence it is ``covering" for the linearized second order elliptic differential operator.

%At the end, we \textcolor{red}{prove} the relation \eqref{eq:idcal}, which can be written as
%\begin{equation}\label{eq:proofApp1}
%	2\,F^{\top}\,\partial_{1}
%	(A_H F)
%	=
%	\partial_{1}\pare{ F^{\top}\,  A_H F} 
%	+
%	F^{\top}\, 
%	\pare{\partial_{1}A_H}\,F
%	.
%\end{equation}
%Since $A_H=A\circ H$ is a real-valued symmetric matrix, we can write $A_H=B^{\top}B$.
%Then 
%$\partial_{1}( F^{\top} A_H F) =2F^{\top}B^{\top}\partial_1(BF)$
%and 
%$F^{\top}\partial_{1}
%(A_H F)
%=F^{\top}B^{\top}\partial_1(BF)+F^{\top}(\partial_1 B^{\top})BF$.
%Hence, by noticing that $2\partial_1(B^{\top})B=\partial_{1}(B^{\top}B)=\partial_{1}A_H$ we have proven \eqref{eq:proofApp1}.

\frenchspacing

\section*{Acknowledgments}
{The research of FC was partially supported  by the AFOSR Grant  FA9550-20-1-0024 and  NSF Grant DMS-21-06255. The research of MSV was partially supported by NSF Grant DMS-22-05912. Data sharing is not applicable to this article as no datasets were generated or analysed during the current study. The authors have no financial or proprietary interests in any material discussed in this article.}

%\bibliographystyle{abbrv}%
%\bibliography{book}

\end{document}